\documentclass[12pt,a4paper]{article}

\usepackage{amsmath,amsfonts}
\usepackage{mathrsfs} 
\usepackage{parskip}
\usepackage[T1]{fontenc}
\usepackage[utf8]{inputenc}
\usepackage{authblk}
\usepackage{enumerate}
\usepackage{graphicx}

\newtheorem{theorem}{Theorem}

\newtheorem{proposition}{Proposition}

\newtheorem{example}{Example}

\let\ds\displaystyle

\def\AA{\mathbb{A}}
\def\BB{\mathbb{B}}

\def\Pb{\mathbf{P}}

\def\Ex{\mathbf{E}}
\def\Pb{\mathbf{P}}

\def\KK{\mathbb{K}}

\def\YY{\mathbb{Y}}  
\def\sgn{{\rm sgn}}
\def\1{\mbox{1\hspace{-.25em}I}}
\begin{document}

\title{On Non Asymptotic Expansion of the MME in the Case of Poisson
  Observations.}

\author[1]{O.V. Chernoyarov}
\author[2]{A.S. Dabye}
\author[3]{F.N. Diop }
\author[4]{Yu.A. Kutoyants}
\affil[1]{\small  National Research University ``MPEI'', Moscow, Russia}
\affil[2]{\small  University Gaston Berger, Saint Louis, Senegal}
\affil[3]{\small University of Thies, Thies, Senegal}
\affil[4]{\small  Le Mans University,  Le Mans,  France}
\affil[1,4]{\small  Tomsk State University, Tomsk, Russia}

\date{}

\maketitle
\begin{abstract}
 The problem of parameter estimation by observations of
 inhomogeneous Poisson processes is considered.   The method of moments
 estimator is studied and its stochastic expansion is obtained. This stochastic
 expansion is then used to obtain the expansion of the moments of this
 estimator and the expansion of the distribution function. The stochastic
 expansion, expansion of the moments and the expansion of distribution
 function are non asymptotic in nature.  Several examples are considered. 
\end{abstract}
\noindent MSC 2000 Classification: 62M05,  62G05, 62G20.

\bigskip
\noindent {\sl Key words}: \textsl{Poisson process, Parameter estimation,
  method of moments, expansion of estimators, expansion of the moments,
  expansion of distribution function, non asymptotic expansions.}

\section{Introduction}

Consider the problem of parameter estimation by $n$
independent observations $X^{(n)}=\left(X_1,\ldots,X_n\right)$. If we suppose
that the density function $f\left(\vartheta ,x\right)$ of the observation
$X_j$ is {\it regular}, i.e., is sufficiently smooth with respect to parameter
$\vartheta $, then the well-know estimators (MLE, Bayesian estimator, method
of moments estimators) are consistent, asymptotically normal  and we have the
convergence of polynomial moments ($p>0$):
\begin{align*}
&\bar \vartheta _n\longrightarrow \vartheta _0, \qquad  \qquad  \sqrt{n}\left(\bar
\vartheta _n- \vartheta _0\right)\Longrightarrow \xi \sim {\cal N}\left(0,{\rm
  D}\left(\vartheta _0\right)^2\right),\\
&\lim_{n\rightarrow \infty }n^{p/2}\Ex_{\vartheta _0}\left| \bar
\vartheta _n- \vartheta _0\right|^p\longrightarrow  {\rm
  D}\left(\vartheta _0\right)^\frac{p}{2}\Ex \left|\zeta \right|^p,\qquad
\zeta \sim {\cal N}\left(0,1\right). 
\end{align*}
Here we denoted $\vartheta _0$ the true value and  ${\rm
  D}\left(\vartheta _0\right)^2$ is the limit variance of the estimator
$\bar\vartheta _n$. These relations can be written as follows
\begin{align}
&\bar \vartheta _n-\vartheta _0=o\left(1\right) ,\nonumber\\
&  \label{1}
\bar
\vartheta _n- \vartheta _0=\varphi _n\xi \left(1+o\left(1\right)\right) ,\\
&\Ex_{\vartheta _0}\left| \bar
\vartheta _n- \vartheta _0\right|^2= \varphi _n^2{\rm
  D}\left(\vartheta _0\right)^2 \left(1+o\left(1\right)\right),\label{2}\\
&\Pb_{\vartheta _0}\left({\rm
  D}\left(\vartheta _0\right)^{-1}\varphi
_n^{-1}\left(\bar \vartheta _n-\vartheta _0\right)<x\right)=
F\left(x\right)+o\left(1\right). \label{3}
\end{align}
Here $\varphi _n=n^{-1/2}$, we take $p=2$ and $F\left(x\right)$ is distribution
function of Gaussian law ${\cal N}\left(0,1\right)$.  Of course, the relation
in \eqref{1} is just a symbolique writing because the limit Gaussian variable
$\xi $ is not defined on the same probability space, we have convergence in
distribution only.

If the volume $n$ of observations is large, then the relations \eqref{1}-\eqref{3}  describe well
the distribution of the  error of estimation. For the
moderate values of $n$ the real distribution of $\bar\vartheta _n-\vartheta
_0$ and of  the moments $\Ex_{\vartheta _0}\left| \bar \vartheta _n- \vartheta
_0\right|^2 $ can be quite far from the given here limit values. The better
approximations for the distribution function and the moments can be obtained
with the help of well-known {\it asymptotic expansion } theory. There are at
least three types of expansions:

{\it Stochastic expansion}
\begin{align}
\label{4}
\bar \vartheta _n- \vartheta _0=\varphi _n\xi _{n,1}+\ldots+\varphi _n^k\xi
_{n,k}+o\left(\varphi _n^k \right),
\end{align}
where $\xi _{n,i}$ are bounded in probability random variables,

{\it Expansion of the moments}
\begin{align}
\label{5}
\Ex_{\vartheta }\left|\bar \vartheta _n- \vartheta _0\right|^2=\varphi _n^2P
_{n,1}+\ldots+\varphi _n^{2k}P _{n,k}+o\left(\varphi _n^k \right),
\end{align}
where $P _{n,i}$ are some bounded real values,

{\it Expansion of the distribution function}
\begin{align}
\label{6}
\Pb_{\vartheta _0}\left(\frac{ \bar \vartheta _n- \vartheta
_0}{{\rm D}\left(\vartheta _0\right)\sqrt{n}}<x\right)=
F\left(x\right)+\varphi _np _1\left(x\right)+\ldots+\varphi
_n^kp_k\left(x\right) +o\left(\varphi _n^k\right),
\end{align}
where   $p_i\left(x\right)$ are some products of polynomials
and $f\left(x\right)$ (density function of ${\cal
  N}\left(0,1\right)$). The value of the  parameter $k$ depends on the
smoothness of the model with respect to unknown paramater.

 We can mention here the works devoted to asymptotic expansions of estimators
 for independent identically distributed observations of the random variables
 \cite{BH90}, \cite{BU82}, \cite{CH72}, \cite{CH73}, \cite{GU75}, \cite{GU76},
 \cite{LM64}, \cite{M67}, \cite{PF82}. In the last book there is an extensive
 list of references. Such asymptotic expansions are widely used in  bootstrap
 too \cite{BSY03}. The difference between these works is in the conditions
 of regularity and in the estimates of the residuals $o\left(\cdot \right)$ in
 \eqref{4}-\eqref{6}.
Note that in the majority of these works the results are {\it asymptotic in
  nature}, i.e., the residuals in expansions \eqref{3}-\eqref{5}   are of the type $o\left(\varphi
_n^k\right)$ or     $ \varphi _n^{k+\delta } o\left(1\right), \delta \in
\left(0,1\right)$, where $o\left(1\right)\rightarrow 
0$ as $n\rightarrow \infty $ and nothing can be said about the term
$o\left(1\right)$ for finite $n$.

The expansions of the errors of estimation by the powers of small parameters
in the case of observations of continuous time stochastic processes were
obtained in the works \cite{BU77} (signal in white Gaussian noise),
\cite{Kut83}, \cite{Kut98} (inhomogeneous Poisson processes),
\cite{Kut84},\cite{Y92}, \cite{Y93}, \cite{Kut94}, (diffusion processes),
\cite{Y01} (martingales with jumps).

One of the goals of this work is to obtain such expansions (stochastic,
moments and distribution function) for the recently introduced class of
estimators: {\it method of moments estimators} (MME) in the case of
observations of inhomogeneous Poisson processes. The method of moments allowed
to have explicit expressions for many models with intensity functions non
linearly depending on the unknown parameters, where the traditional MLE have
no explicit expressions and therefore MME provides essential gain in the
calculation of the consistent and asymptotically normal estimators. Note that
the MME are used in One-step MLE construction to obtain the asymptotically
efficient estimators too \cite{DGK18}.

 It is known that the publications
related with the asymptotic expansions \eqref{4}-\eqref{6} in statistics are
technically quite cumbersome (see, e.g., \cite{PF82} and references
therein). In the considered in this work case the exposition is essentially
simplified because the random variables $\xi _{n,i}$ in \eqref{4} have the
form $\xi _{n,i}=\eta _n^i$ with the same $\eta _n$. Another advantage of this
work with respect to traditional asymptotic expansions is the using in
obtaining all expansions of the {\it method of good sets}, which in reality
allows to obtain {\it non asymptotic expansions}. This method was developed by
Burnashev in \cite{BU77} and \cite{BU82}. His approach allows to have
expansions \eqref{4} and \eqref{5} {\it non asymptotic in nature}. This means,
that the residuals in \eqref{4} and \eqref{5} for finite $n$ can be estimated
from above and from below with known rates and constants. Remark that this
method was already applied to obtain non asymptotic expansions in the works
\cite{Kut83}, \cite{Kut84}, \cite{Kut94}, \cite{Kut98}.

\section{Method of Moments Estimator}

Let us consider the following problem of parameter estimation. Suppose that we
have $n$ independent observations of inhomogeneous Poisson processes
$X^{(n)}=\left(X_{1},\ldots,X_{n}\right)$, where $X_{j}\!=\!\bigl\{
X_{j}(t),\, t\in {\cal T} \bigr\}$ with intensity function
$\lambda\left(\vartheta,t\right),t\in {\cal T}$. Here $ {\cal T}\subset {\cal
  R}$ is some interval. It can be ${\cal T}=\left[0,\tau \right]$, ${\cal
  T}=[0,+\infty ) $, $ {\cal T}= {\cal R}$ or any other interval on the
  line. As usual, $X_j\left(\cdot \right)$ are counting processes (see details
  in \cite{Kut98}).  We have to estimate the parameter $\vartheta \in\Theta
  =\left(\alpha ,\beta \right)$ by the observations $X^{\left(n\right)}$ and
  to describe the properties of estimators in the asymptotic of large samples
  ($n\rightarrow \infty $).

For the study we take the method of moments estimator (MME) defined as
follows.  Introduce the functions $g(t), t\in {\cal T}$ and the functions
\[
m(\vartheta)=\int_{\cal T} g(t)\lambda(\vartheta,t){\rm d}t,\quad
\vartheta\in\Theta\quad 
\mbox{and}\quad \bar{m}_n=\frac{1}{n}\sum_{j=1}^n\int_{\cal T}g(t){\rm d}X_j(t).
\]
We suppose that we have such  functions $\lambda \left(\vartheta
,t\right),\vartheta \in\Theta ,t\in {\cal T},$
and $g\left(t\right)$ that the function $m\left(\vartheta \right),\alpha \leq
\vartheta \leq \beta $ is monotone. Without loss of generality we assume that
it is monotone increasing.

Let us introduce the following notations:
\begin{align*}
{\cal M}=\left\{m\left(\vartheta \right):\; \vartheta \in \left[\alpha ,\beta
  \right]\right\}= \left[m\left(\alpha \right),m\left(\beta \right)\right]
\end{align*}
For $y\in{\cal M}$ we write the solution of the equation $m(\vartheta)=y$  as 
$\vartheta=m^{-1}(y)=G(y)$. Therefore the function $G(y)$ is inverse 
for $m(\vartheta)$ and $G(m(\vartheta))=\vartheta$.

We define the {\it method of moments estimator} (MME) $\check{\vartheta}_n$
by the following equation 
\[
\check{\vartheta}_n=\arg
\inf_{\vartheta\in\Theta}\left(m(\vartheta)-\bar{m}_n\right)^2.
\]
This estimator admits the representation
\begin{align}
\label{2-1}
\check{\vartheta}_n=\alpha\1_{\left\{\bar{m}_n\leq m(\alpha)\right\}}
+\bar{\vartheta}_n\1_{\left\{m(\alpha)<\bar{m}_n<m(\beta)\right\}}+\beta\1_{\left\{\bar{m}_n\geq m(\beta)\right\}},
\end{align}
where $\bar{\vartheta}_n$ is a solution of the equation $m(\bar{\vartheta}_n)=\bar{m}_n$
when $\bar{m}_n\in\left(m(\alpha),m(\beta)\right)$. The value
$\check{\vartheta}_n=\alpha$ (respectively 
$\check{\vartheta}_n=\beta$) corresponds to $m\leq m(\alpha)$ closest to
$\bar{m}_n$ (respectively $m\geq m(\beta)$ closest to $\bar{m}_n$). Recall
that with probability $p_n=\exp\left(-n\Lambda \left({\cal T}\right)\right)$
we have no events (jumps) in all observations. Here 
$$
\Lambda \left({\cal
  T}\right)=\int_{{\cal T}}^{}\lambda \left(\vartheta ,t\right){\rm d}t.
 $$   

This estimator was recently studied in \cite{DGK18}. It was shown that as
$n\rightarrow \infty $ and under mild regularity conditions this estimator is
consistent
\begin{align}
\label{2-2}
\check{\vartheta}_n\longrightarrow \vartheta _0
\end{align}
(here and in the sequel $\vartheta _0$ denotes the true value), asymptotically
normal
\begin{align}
\label{2-3}
\sqrt{n}\left(\check{\vartheta}_n-\vartheta _0\right)\Longrightarrow {\cal
  N}\left(0, {\rm D}\left(\vartheta _0\right)^2\right)
\end{align}
or, for any $x\in{\cal R}$
\begin{align*}
\Pb_{\vartheta _0}\left( {\rm D}\left(\vartheta
_0\right)^{-1}\sqrt{n}\left(\check{\vartheta}_n-\vartheta 
_0\right)<x  \right)\longrightarrow  F\left(x\right)=\frac{1}{\sqrt{2\pi
}}\int_{-\infty }^{x}e^{-\frac{y^2}{2}}{\rm d}y . 
\end{align*}
We have as well the convergence of all polynomial moments: for any $p>0$
\begin{align}
\label{2-4}
n^{\frac{p}{2}}\Ex_{\vartheta _0}\left| \check{\vartheta}_n-\vartheta _0
\right|^p\longrightarrow  \Ex\left|\zeta \right|^p\; {\rm D}\left(\vartheta
_0\right)^{\frac{p}{2}}, \qquad \zeta \sim {\cal N}\left(0,1\right). 
\end{align}
The advantage of the MME can be illustrated with the help of  the following examples.

\begin{example} {\rm Suppose that we have an inhomogeneous Poisson process
    $X^T=\left(X_t,0\leq t\leq T\right) $ with $\tau $-periodic intensiy
    function $\lambda \left(\vartheta,t \right),0\leq t\leq T$, i.e., $\lambda
    \left(\vartheta,t+k\tau  \right)=\lambda \left(\vartheta,t  \right) $ for
    any $k=1,2,\ldots$. Suppose that
    $T=n\tau $ and introduce the independent Poisson processes
\begin{align*}
X_j\left(t\right)=X_{\left(j-1\right)\tau +t}-X_{\left(j-1\right)\tau},\quad
t\in {\cal T}=\left[ 0, \tau\right] ,\quad j=1,\ldots,n. 
\end{align*}
Therefore we obtain the observations
$X^{\left(n\right)}=\left(X_1,\ldots,X_n\right)$ with intensity function
$\lambda \left(\vartheta ,t\right)$ and can
construct the MME $\check{\vartheta}_n $.

 Suppose that we observe a $\tau $-periodic Poisson
signal of intensity function 
$S\left(\vartheta ,t\right)=\vartheta h\left(t\right)$ (amplitude modulation)
in the presence of the Poisson noise
of intensity $ \lambda _0>0$, i.e., we have
\begin{align}
\label{2-5}
\lambda \left(\vartheta ,t\right)=\vartheta h\left(t\right)+\lambda _0,\qquad 
0\leq t\leq \tau . 
\end{align}
Here $h\left(\cdot \right)$ is $\tau $-periodic known positive function and
$\lambda _0>0$. Remind that the MLE $\hat\vartheta _n$ for this model has no
explicit representation and is given as solution of the following equation
\begin{align*}
\sum_{j=1}^{n}\int_{0}^{\tau }\frac{h\left(t\right)}{\hat\vartheta
  _nh\left(t\right)+\lambda _0 }{\rm d}X_j\left(t\right)=n\int_{0}^{\tau
}h\left(t\right){\rm d}t. 
\end{align*}
Despite the numerical difficulties of its calculation there is as well the
problem of definition of this stochastic integral. The calculation of the MME
has no such difficulties and for the wide class of functions $g\left(\cdot
\right)$ (say, positive)  we have
\begin{align*}
m\left(\vartheta \right)=\vartheta  H_g+  \lambda _0 G,\qquad H_g=\int_{0}^{\tau
}g\left(t\right)h\left(t\right){\rm d}t,\quad G=\int_{0}^{\tau
}g\left(t\right){\rm d}t
\end{align*}
and $ \vartheta =H_g^{-1}\left[m\left(\vartheta \right)-\lambda _0
  G\right]$. Hence $ G\left(y\right)=H_g^{-1}\left[y-\lambda _0
  G\right]  $. 
The MME
\begin{align}
\label{2-7}
\check{\vartheta}_n=\left(\int_{0}^{\tau }g\left(t\right)h\left(t\right){\rm
  d}t\right)^{-1}\left[ \frac{1}{n}\sum_{j=1}^{n}\int_{0}^{\tau
  }g\left(t\right){\rm d}X_j\left(t\right)   -\lambda _0 \int_{0}^{\tau 
  }g\left(t\right){\rm d}t \right]
\end{align}
This estimator has the mentioned above properties \eqref{2-2}-\eqref{2-4} (see
\cite{DGK18}).

 }
\end{example}
\begin{example}. {\rm Suppose that the intensity function of the observed
    Poisson processes $X^{\left(n\right)}=\left(X_1,\ldots,X_n\right)$,
    $X_j=\left(X_j\left(t\right), t\in {\cal T}=[0,+\infty ) \right)$ is
\begin{align}
\label{2-8}
\lambda \left(\vartheta ,t\right)=f\left(t\right) \;e^{-\vartheta h\left(t\right)}+q\left(t\right),\quad t\geq 0.
\end{align}
Here $f\left(\cdot \right)>0, h\left(\cdot \right)\geq 0$ and $q\left(\cdot
\right)\geq 0 $
are known 
functions, $h\left(0\right)=0,h\left(\infty \right)=\infty $ and $h\left(\cdot
\right)$ has  continuous derivative $h'\left(\cdot \right)$. For example,
$f\left(t\right)=1  +t^4, h\left(t\right)=t^2,q\left(t\right)=1$.
  We have to estimate $\vartheta \in\left(\alpha ,\beta \right)$, $0<\alpha <\beta <\infty $. Recall
 that the MLE has no explicit expression. Let us put
 $g\left(t\right)=f\left(t\right)^{-1} h'\left(t\right)$, then
 we have 
\begin{align*}
m\left(\vartheta \right)=\int_{0}^{\infty }g\left(t\right)\lambda \left(\vartheta
,t\right){\rm d}t=\frac{1}{\vartheta }+R, \quad R=\int_{0}^{\infty }  \frac{h'\left(t\right)q\left(t\right)}{f\left(t\right)}{\rm d}t,\quad   
\end{align*}
hence  $\vartheta= \left(m\left(\vartheta \right)-R\right)^{-1}  $   and the MME 
\begin{align*}
\check{\vartheta}_n=\left(\frac{1 }{n}\sum_{j=1}^{n}\int_{0}^{\infty
}\frac{h'\left(t\right)}{f\left(t\right)}{\rm d}X_j\left(t\right)-R
\right)^{-1}.
\end{align*}
If we suppose that the corresponding integrals are finite then once more this
estimator has the mentioned above properties \eqref{2-2}-\eqref{2-4} (see
\cite{DGK18}).

The case $h\left(0\right)\not=0$ can be treated by a similar way, but to
define the MME we have to solve the equation 
\begin{align*}
\vartheta ^{-1}e^{-\vartheta h\left(0\right)}=y+R,\qquad \vartheta \in\Theta .
\end{align*} 
Of course, all derivatives of the solution $\vartheta =G\left(y\right)$ can be
calculated without problems. 
}

\end{example}
\begin{example} {\rm Consider Poisson processes
    $X^{\left(n\right)}=\left(X_1,\ldots,X_n\right)$ with ``Gaussian'' intensity
    function
\begin{align}
\label{2-9}
\lambda \left(\vartheta ,t\right)=a\;e^{-\frac{\left(t-b\right)^2}{2\vartheta
    ^2}},\quad t\in {\cal R}.
\end{align}
Here $a>0, b $ are supposed to be  known and we have to estimate $\vartheta \in
\left(\alpha ,\beta \right)$, $\alpha >0$. Let us take
$g_1\left(t\right)=\left(t-b \right)^2$. Then
\begin{align*}
m_1\left(\vartheta \right)=\int_{\cal R}^{}\left(t-b \right)^2\lambda
\left(\vartheta ,t\right){\rm d}t =\vartheta ^3a\sqrt{2\pi },\qquad
 \vartheta =\left(\frac{m\left(\vartheta \right)}{a\sqrt{2\pi }}\right)^{1/3}
\end{align*} 
and 
\begin{align*}
\check{\vartheta}_n=\left(\frac{1 }{an\sqrt{2\pi }}\sum_{j=1}^{n}\int_{\cal
  R}^{}\left(t-b \right)^2 {\rm d}X_j\left(t\right) \right)^{1/3}.
\end{align*}
Another possibility is to take $g_2\left(t\right)=\left|t-b \right|$. Then
we obtain $m_2\left(\vartheta \right)=2a\vartheta ^2$ and
\begin{align*}
\check{\vartheta}_n=\left(\frac{1 }{2an}\sum_{j=1}^{n}\int_{\cal R}^{}\left|t-b  \right| {\rm
  d}X_j\left(t\right)\right)^{1/2}.
\end{align*}
The both estimators are consistent and asymptotically normal.
}
\end{example}

\section{Stochastic expansion}

The properties \eqref{2-3}, \eqref{2-4} of the MME $\check{\vartheta}_n$ can
be written as \eqref{1}, \eqref{2} and our goal is to obtain the expansions like
\eqref{4} for these estimators too.

Recall, we have $n$ independent observations of inhomogeneous Poisson
processes $X^{\left(n\right)}$ of the same  intensity function $\lambda
\left(\vartheta ,t\right), t\in {\cal T} $, where $\vartheta \in\Theta =\left(\alpha
,\beta \right) $.  To estimate $\vartheta $ we use the MME
$\check{\vartheta}_n $ defined in \eqref{2-1} with some  function
$g\left(\cdot \right)$.

Introduce the notation: 
\begin{align*}
\psi_l\left(\vartheta _0\right)&=\frac{G^{(l)}(m(\vartheta_0))}{l!}, l=1,\ldots,k,\quad \pi
_n\left(t\right)= \sum_{j=1}^{n}X_j\left(t\right)-n\int_{s<t}^{}\lambda
\left(\vartheta _0,s\right){\rm d}s
,\\
\eta _n\left(\vartheta _0\right)&=\sqrt{n}\left(\bar m_n-m\left(\vartheta
_0\right)\right)=\frac{1}{\sqrt{n}}\int_{ {\cal T}}^{
}g\left(t\right)\,{\rm     d}\pi _n\left(t\right) .
\end{align*}
Remind that sufficient condition for  the asymptotic normality of the
stochastic integrals $\eta _n\left(\vartheta \right),\vartheta \in\Theta  $ is
\begin{align*}
\sup_{\vartheta \in\Theta }\int_{ {\cal T}}^{ }g\left(t\right)^2\,\lambda \left(\vartheta ,t\right){\rm d}t<\infty .
\end{align*}
Below we introduce more strong condition $\mathcal{L}_3 $.

 In this work we suppose that the functions  $\lambda
\left(\vartheta ,t\right)$ and $g\left(t\right)$ are such that the following
conditions $\mathcal{L}$  are fulfilled.
\begin{description}
\item[$\mathcal{L}_1.$] {\it The function $m(\vartheta),\,\vartheta\in\Theta$
has $k+2$ continuous bounded derivatives.

\item[$\mathcal{L}_2.$] The function
 $m(\vartheta),\,\vartheta\in\Theta$ is monotone and
\[
\inf_{\vartheta\in\Theta}|\dot{m}(\vartheta)|>0.
\] 
\item[$\mathcal{L}_3.$] The function $g\left(\cdot \right)$ is such  that for
  any $m>0$ }
\begin{align}
\label{3-0}
\sup_{\vartheta \in\Theta }\int_{ {\cal T}}^{ }\left|g\left(t\right)\right|^m\,\lambda \left(\vartheta ,t\right){\rm d}t<\infty .
\end{align}
\end{description}
Without loss of generality, we suppose that the function $m(\vartheta ),
\vartheta \in\Theta$ is increasing.  We have the following first result
concerning the stochastic expansion of the MME.
\begin{theorem}
\label{T3-1}
Let the conditions $\mathcal{L}$ be fulfilled, then for any  $k=1,2,\ldots $ there exist the random
variables $r_{n,k}$, $\phi_{n,k} $ and the set $\mathbb{B}$ such that
 the MME $\check{\vartheta}_n$ admits the representation
{\rm
\begin{equation}
\label{3-1}
\check{\vartheta}_n=\vartheta_0+
\left\{\sum_{l=1}^k\psi_l\left(\vartheta
_0\right)\,\eta_n^l\;n^{-\frac{l}{2}}+r_{n} n^{-\frac{k}{2}-\frac{1}{4}} 
\right\}\1_{\{\mathbb{B}\}}+\phi_n\1_{\{\mathbb{B}^c\}},
\end{equation}
}where $\eta _n=\eta _n\left(\vartheta _0\right), \vert r_n \vert\leq  1$, $\phi_n\in \left(\alpha -\vartheta
_0,\beta -\vartheta _0\right)$  and  for any $Q>0 $ and any compact
$\KK\subset\Theta $ there exists a
constant $C=C\left(Q,\KK\right)>0$ such that
\begin{align}
\label{3-2}
\sup_{\vartheta \in \KK}\Pb_{\vartheta_0}\left(\mathbb{B}^c\right)\leq \frac{C}{n^Q}.
\end{align}
\end{theorem}
{\bf Proof.}  
The proof of this theorem is based on the approach of {\it good sets}. 
Introduce the first {\it good set}
\[
\mathbb{B}_1=\left\{
\inf_{|\vartheta-\vartheta_0|<\delta}
\vert m(\vartheta)-\bar{m}_n \vert <
\inf_{|\vartheta-\vartheta_0|\geq\delta}
\vert m(\vartheta)-\bar{m}_n \vert 
\right\},
\]
where $\delta >0$ is some small number satisfying the condition $ \alpha
+\delta <\vartheta _0<\beta -\delta $. 
Then the MME \eqref{2-1} on the set $\mathbb{B}_1$ satisfies the relations
\[
m(\check{\vartheta}_n)=\bar{m}_n \qquad 
\check{\vartheta}_n=G(\bar{m}_n)=G\left( m(\vartheta_0)+\varepsilon\eta_n\right).
\]  
Here $\varepsilon=n^{-1/2}$  and  $\bar{m}_n= m(\vartheta_0)+n^{-1/2}\eta_n $.

The Taylor expansion of the function $G(\cdot)$ on the set $\mathbb{B}_1$ yields
\begin{eqnarray*}
\check{\vartheta}_n&=&G\left( m(\vartheta_0)\right)+
\sum_{l=1}^k\frac{G^{(l)}\left(
  m(\vartheta_0)\right)}{l!}\,\eta_n^l\,\varepsilon^l+ \frac{G^{(k+1)}(
  \tilde{m}_n)}{(k+1)
}{\eta_n^{k+1}}\,\varepsilon^{k+1}\\ &=&\vartheta_0+\sum_{l=1}^k\frac{G^{(l)}\left(
  m(\vartheta_0)\right)}{l!}\eta_n^l\left(\frac{1}{\sqrt{n}}\right)^l
+\frac{G^{(k+1)}(
\tilde{m}_n)}{ (k+1)!  }\frac{\eta_n^{k+1}}{n^{\frac{1}{4}}}\left(\frac{1}{\sqrt{n}}\right)^{k+\frac{1}{2}}\\
 &=&\vartheta_0+\sum_{l=1}^k\psi
_l(\vartheta_0)\;{\eta_n^l}\;n^{-\frac{l}{2} } +r_{n,k}\;n^{-\frac{k}{2}-\frac{1}{4}},
\end{eqnarray*}
where $\tilde{m}_n\in\left(m(\vartheta_0-\delta),m(\vartheta_0+\delta)\right)$
and we denoted  
$$
r_{n,k}= \ds\frac{ 
  G^{(k+1)}( \tilde{m}_n)\;\eta_n^{k+1}}{(k+1)!\;n^{\frac{1}{4}}}. 
$$

Recall that the derivatives $G'(y)$, $G''(y)$, $G'''(y)$ of the inverse 
function $G(y)$ can be calculated using the equality $G(m(\vartheta))=\vartheta$
as follows
\begin{multline*}
G'(m(\vartheta))\dot{m}(\vartheta)=1,\qquad 
G'(m(\vartheta))=\frac{1}{\dot{m}(\vartheta)},\qquad
G'(y)=\frac{1}{\dot{m}(G(y))}, \\
G''(y)=-\frac{\ddot{m}(G(y))}{\dot{m}(G(y))^3},\qquad
G'''(y)=\frac{3\ddot{m}(G(y))^2-\dot{m}(G(y))
\stackrel{...}{m}(G(y))}{\dot{m}(G(y))^5}.
\end{multline*}
Here dot means derivation w.r.t. $\vartheta $.  The other derivatives can be calculated by the same rules. As it follows from
conditions ${\cal L}$ all derivatives up to $G^{(k+1)}( \cdot ) $ are bounded
and we can write on the set $\BB_1$
(with probability 1)
\begin{align*}
\left|\left.\frac{\partial ^{k+1}G(y)}
{\partial y^{k+1}}\right|_{y=\tilde{m}_n}\right|
\leq
\sup_{m\left(\vartheta _0 -\delta \right)\leq y\leq m\left(\vartheta _0
  +\delta \right)}  \left|\frac{\partial ^{k+1}G(y)}
{\partial y^{k+1}}\right|   =C_{k+1},
\end{align*}
where the constant $C_{k+1}=C_{k+1}\left(\vartheta _0,\delta \right)>0 $.

 Introduce the second \textsl{good set}
\[\mathbb{B}_2=
\left\{\,\vert \eta_n\vert^{k+1}<n^{1/4}\frac{(k+1)!}{C_{k+1}}
\right\}.
\]
 On the set   $\BB=\BB_1\cap\BB_2$ we have the estimate $\left|r_{n,k}\right|\leq 1 $.

Therefore we obtain the stochastic expansion \eqref{3-1} on the set
$\mathbb{B}$.   The probability of the
complement is estimated as follows
\[
\Pb_{\vartheta_0}\left(\mathbb{B}^c\right)\leq \Pb_{\vartheta_0}\left(\mathbb{B}_1^c\right)+
\Pb_{\vartheta_0}\left(\mathbb{B}_2^c\right).
\]

To estimate $\Pb_{\vartheta_0}\left(\mathbb{B}_1^c\right)$ 
we remark that 
\begin{align*}
\Pb_{\vartheta_0}\left(\mathbb{B}_1^c\right)&=
\Pb_{\vartheta_0}\left(\vert
\check{\vartheta}_n-\vartheta_0\vert\geq \delta\right)\\
&=
\Pb_{\vartheta_0}\left(
\inf_{|\vartheta-\vartheta_0|<\delta}\vert m(\vartheta)-
\bar{m}_n\vert\geq
\inf_{|\vartheta-\vartheta_0|\geq\delta}\vert m(\vartheta)-
\bar{m}_n\vert\right)\\
&=
\Pb_{\vartheta_0}\left(
\inf_{|\vartheta-\vartheta_0|<\delta}\vert m(\vartheta)-m(\vartheta_0)+m(\vartheta_0)
-\bar{m}_n\vert  \right.\\
& \qquad \qquad \qquad \left.\geq
\inf_{|\vartheta-\vartheta_0|\geq\delta}\vert m(\vartheta)-m(\vartheta_0)+m(\vartheta_0)-
\bar{m}_n\vert\right)\\
&\leq
\Pb_{\vartheta_0}\left(2\vert \bar{m}_n-m(\vartheta_0)\vert \geq 
\inf_{|\vartheta-\vartheta_0|\geq\delta} \vert
m(\vartheta)-m(\vartheta_0)\vert\right)\\
&\leq
\Pb_{\vartheta_0}\left(2\vert \bar{m}_n-m(\vartheta_0)\vert \geq \rho \left(\delta \right)\right),
\end{align*}
where we denoted $\rho \left(\delta \right)= \inf_{|\vartheta-\vartheta_0|\geq\delta} \vert
m(\vartheta)-m(\vartheta_0)\vert$. Note that 
 for any $\delta >0$  by condition $\mathcal{L}_2$  we  can write 
\begin{equation}
\label{3-4}
\rho(\delta )=\inf_{|\vartheta-\vartheta_0|>\delta }
\vert \dot{m}(\tilde{\vartheta})\vert |\vartheta-\vartheta_0|
\geq\kappa\delta >0.
\end{equation}
Here 
$\kappa=\inf_{\vartheta\in\Theta}\vert\dot{m}(\vartheta)\vert$.
Therefore  we have 
\begin{align*}
\Pb_{\vartheta_0}\left(\mathbb{B}_1^c\right)&\leq
\Pb_{\vartheta_0}   \left(\left|\frac{2}{n}\sum_{j=1}^{n}\int_{{\cal T}}^{}
g\left(t\right)\left[{\rm d}X_j\left(t\right)-\lambda \left(\vartheta
  _0,t\right){\rm d}t\right]\right| \geq \kappa \delta \right) \\
& \leq  \Pb_{\vartheta_0}   \left(\frac{2}{\sqrt{n}}\left|\int_{{\cal T}}^{ }
g\left(t\right){\rm d}\pi _n\left(t\right)\right| \geq \kappa \delta\sqrt{n}
\right),
\end{align*}
where 
$$
{\rm d}\pi _n\left(t\right)={\rm d}   Y_n\left(t\right)-n\lambda \left(\vartheta
  _0,t\right){\rm d}t,\qquad \quad Y_n\left(t\right)=\sum_{j=1}^{n}  X_j\left(t\right).
$$
Remark that $Y_n\left(t\right), t\in {\cal T} $ is inhomogeneous Poisson
process with intensity function $n\lambda \left(\vartheta _0,t\right), t\in
{\cal T} $.

Further, by Markov inequality we can write for any integer $q\geq 1$
\begin{align}
\label{3-5}
&\Pb_{\vartheta_0}   \left({\frac{2}{\sqrt{n}}}\left|\int_{{\cal T}}^{}
g\left(t\right){\rm d}\pi _n\left(t\right) \right|\geq \kappa \delta \sqrt{{n}}
\right)\nonumber\\
&\qquad\qquad\leq \left(\kappa \delta \sqrt{{n}}\right)^{-2q}\Ex_{\vartheta
  _0}\left|\int_{{\cal T}}^{} 
g\left(t\right){\rm d}\pi _n\left(t\right)
\right|^{2q}\nonumber\\
&\qquad\qquad\leq C_14^q\left(\kappa \delta
\sqrt{{n}}\right)^{-2q}n^{1-q}\int_{{\cal T}}^{} 
\left|g\left(t\right)\right|^{2q} \lambda \left(\vartheta _0,t\right)   {\rm
  d}t\nonumber\\
&\qquad \qquad \qquad + C_24^q\left(\kappa \delta
\sqrt{{n}}\right)^{-2q}\left(\int_{{\cal T}}^{} 
\left|g\left(t\right)\right|^{2} \lambda \left(\vartheta _0,t\right)   {\rm
  d}t\right)^q  \leq \frac{C}{n^q}.
\end{align}
Here we used the property of stochastic integral
\begin{align}
\label{3-6}
\Ex \left|\int_{}^{}f\left(t\right){\rm d}\pi \left(t\right)\right|^{2q}\leq
C_1 \int_{}^{}\left|f\left(t\right)\right|^{2q}\lambda \left(t\right){\rm d}t
+C_2 \left(\int_{}^{}\left|f\left(t\right)\right|^{2}\lambda
\left(t\right){\rm d}t \right)^q 
\end{align} 
with obvious notation. The proof can be found, for example, in \cite{Kut98},
Lemma 1.2.

Further, for the second probability we have
\begin{align}
\label{3-7}
\Pb_{\vartheta_0}\left(\mathbb{B}_2^c\right)&=\Pb_{\vartheta_0}\left(\left|\eta
_n\right|\geq cn^{\frac{1}{4k+4}}\right)\leq
c^{-2q}n^{-\frac{2q}{4k+4}}\Ex_{\vartheta _0}\left|\eta
_n\right|^{2q}\nonumber\\ 
&\leq
C_1c^{-2q}n^{-\frac{2q}{4k+4}}n^{1-q}\int_{{\cal
    T}}^{}\left|g\left(t\right)\right|^{2q}\lambda \left(\vartheta
_0,t\right){\rm d}t\nonumber\\
 &\qquad \quad
+C_2c^{-2q}n^{-\frac{2q}{4k+4}}\left(\int_{{\cal
    T}}^{}\left|g\left(t\right)\right|^{2}\lambda \left(\vartheta
_0,t\right){\rm d}t\right)^q\leq \frac{C}{n^{\frac{q}{2k+2}}}
\end{align}
with the corresponding constant $C>0$.  From the   estimates \eqref{3-5} and
\eqref{3-7} we obtain
\begin{align}
\label{3-8}
\Pb_{\vartheta_0}\left(\mathbb{B}^c\right)&\leq
\Pb_{\vartheta_0}\left(\mathbb{B}_1^c\right)+\Pb_{\vartheta_0}\left(\mathbb{B}_2^c\right)
\leq \frac{C}{n^{\frac{q}{2k+2}}}+ \frac{C}{n^{q}}\leq \frac{C_*}{n^{Q}},
\end{align}
where for a given $Q$ we put $q=2Q\left(k+1\right)$. 
 Note that having functions $g\left(t \right)$ and $\lambda
\left(\vartheta ,t\right)$ all mentioned constants can be calculated or 
estimated from above. 

Remark, that the random variable $\phi _{n,k}$ can take any value on the intervals
$\left[\alpha -\vartheta _0,-\delta \right]$ and $\left[\beta -\vartheta _0,\delta
  ,\right]$. In any case we have the estimate $\left|\phi_n\right|\leq \beta
-\alpha $. It is important to mention that the representation \eqref{3-1} is valid
for all $n$. 

 \section{Expansion of the moments}

The stochastic expansion \eqref{3-1} we write as
\begin{align*}
\check\vartheta _n=\vartheta _0+\Psi _n\1_{\left\{\BB\right\}}+\phi_{n,k} \1_{\left\{\BB^c\right\}},
\end{align*}
where $\Psi _n$ can be written  as follows
\begin{align*}
\Psi _n=\Psi _{0,n}+r_{n,k}\,n^{-\frac{2k+1}{4}}.
\end{align*}
Here of course
\begin{align*}
\Psi _{0,n}=\sum_{l=1}^k\psi_l\left(\vartheta
_0\right)\,\eta_n^l\;n^{-\frac{l}{2}}.
\end{align*}

This presentation allows to write the expansion of the mean of the MME
\begin{align*}
\Ex_{\vartheta _0}\check\vartheta _n=\vartheta _0+\sum_{l=1}^k\psi_l\left(\vartheta
_0\right)\,\Ex_{\vartheta _0}\eta_n^l\;n^{-\frac{l}{2}}+ O\left(n^{-\frac{2k+1}{4}}\right).
\end{align*}
The first terms are
\begin{align*}
\Ex_{\vartheta _0}\check\vartheta _n&=\vartheta _0+\psi _2\left(\vartheta
_0\right)\Ex_{\vartheta _0}\eta_n^2n^{-1} +\psi _3\left(\vartheta
_0\right)\Ex_{\vartheta _0}\eta_n^3n^{-3/2} + O\left(n^{-\frac{7}{4}}\right)\\
&=\vartheta _0-\frac{\ddot m\left(\vartheta _0\right)}{2\dot m\left(\vartheta _0\right)^3} \int_{0}^{T}g\left(t\right)^2\lambda \left(\vartheta _0,t\right){\rm
  d}t\; \frac{1}{n}   + O\left(n^{-\frac{7}{4}}\right).
\end{align*}
We see that if   $k=3$ then the terms with $\psi _1\left(\vartheta _0\right)$ and $\psi
_3\left(\vartheta _0\right)$ are absent in this expansion. These
relations can be proved, but we will speak about more interesting problem of
the expansion of the moments of the error of estimation.

\begin{theorem}
\label{T4-1} Let conditions ${\cal L}$ be fulfilled. Then for any $p>1$ there exists a
constant $C^*>0$  such that 
\begin{align}
\label{4-1}
\Bigl|\Ex_{\vartheta _0}\left|\check\vartheta _n-\vartheta
_0\right|^p-\Ex_{\vartheta _0}\left|\Psi _{0,n}\right|^p\Bigr|\leq C^*\,n^{-\frac{2k+1}{4}}.
\end{align}
\end{theorem}
{\bf Proof.} 
Without loss of generality  we put $r_{n,k}=0$ on the set  $\BB^c$.
Then for any $p\geq 1$  we have  
\begin{align*}
\Ex_{\vartheta _0}\left|\check\vartheta _n-\vartheta
_0\right|^p&=\Ex_{\vartheta _0}\left(\left|\Psi
_n\right|^p\1_{\left\{\BB\right\}} \right)+ \Ex_{\vartheta
  _0}\left(\left|\phi _{n,k}\right|^p\1_{\left\{\BB^c\right\}} \right)\\
&=\Ex_{\vartheta _0}\left|\Psi
_n\right|^p+ \Ex_{\vartheta
  _0}\left(     \left|\phi _{n,k}\right|^p- \left|\Psi
_n\right|^p  \right)\1_{\left\{\BB^c\right\}} .
\end{align*}
Hence
\begin{align*}
\Ex_{\vartheta _0}\left|\check\vartheta _n-\vartheta
_0\right|^p&\leq \Ex_{\vartheta _0}\left|\Psi
_n\right|^p+ \Ex_{\vartheta
  _0}\left(     \left|\phi _{n,k}\right|^p\1_{\left\{\BB^c\right\}}\right)
\end{align*}
and
\begin{align*}
\Ex_{\vartheta _0}\left|\check\vartheta _n-\vartheta
_0\right|^p&\geq \Ex_{\vartheta _0}\left|\Psi
_n\right|^p- \Ex_{\vartheta
  _0} \left( \left|\Psi
_n\right|^p  \1_{\left\{\BB^c\right\}}\right).
\end{align*}

We have
\begin{align*}
\Ex_{\vartheta
  _0}\left(     \left|\phi _{n,k}\right|^p\1_{\left\{\BB^c\right\}}\right)\leq
\frac{C_*\left(\beta -\alpha \right)^p}{n^Q},
\end{align*}
where we used the estimate \eqref{3-8}. 

Note that for any $p>0$ there exist constants $A>0$ and $B>0$ that 
\begin{align*}
\Ex_{\vartheta _0}\left|\eta 
_n\right|^{2p}<A,\qquad \qquad \Ex_{\vartheta _0}\left|\Psi
_n\right|^{2p}<B.
\end{align*}

  By Cauchy-Schwarz inequality we have the estimate
\begin{align*}
\Ex_{\vartheta _0} \left( \left|\Psi _n\right|^p
\1_{\left\{\BB^c\right\}}\right)\leq \left[\Ex_{\vartheta _0} \left|\Psi
  _n\right|^{2p } \Pb_{\vartheta
    _0}\left(\BB^c\right)\right]^{\frac{1}{2}}\leq \frac{\sqrt{C_*B}}{n^{Q/2}}.
\end{align*}
These two estimates allow us to write
\begin{equation}
\label{4-2}
 -\frac{\sqrt{C_*B}}{n^{Q/2}} \leq \Ex_{\vartheta
   _0}\left|\check\vartheta _n-\vartheta_0\right|^p-\Ex_{\vartheta _0}
 \left|\Psi _n\right|^p \leq \frac{C_*\left(\beta -\alpha \right)^p}{n^Q}.
\end{equation}
For $p>1$ and small $x$ we have 
\begin{align*}
\left|a+x\right|^p=\left|a\right|^p+p\;\sgn\left(a+\tilde x\right) \left|a+\tilde
x\right|^{p-1} x,\qquad \left|\tilde x\right|\leq \left|x\right|.
\end{align*}
Using this relation we can write
\begin{align*}
\Bigl|\Ex_{\vartheta _0} \left|\Psi _n\right|^p- \Ex_{\vartheta _0} \left|\Psi
_{0,n}\right|^p \Bigr| & \leq p\Ex_{\vartheta _0}\left|\Psi
_{0,n}+\tilde r_n\,n^{-\frac{2k+1}{4}} \right|^{p-1}\,n^{-\frac{2k+1}{4}}\leq C\,n^{-\frac{2k+1}{4}}.
\end{align*}
Hence
\begin{align*}
 \Ex_{\vartheta
   _0}\left|\check\vartheta _n-\vartheta_0\right|^p-\Ex_{\vartheta _0}
 \left|\Psi _{0,n}\right|^p \leq \frac{C_*\left(\beta -\alpha \right)^p}{n^Q}+C\,n^{-\frac{2k+1}{4}}
\end{align*}
and
\begin{align*}
  \Ex_{\vartheta
   _0}\left|\check\vartheta _n-\vartheta_0\right|^p-\Ex_{\vartheta _0}
 \left|\Psi _{0,n}\right|^p \geq -C\,n^{-\frac{2k+1}{4}}-\frac{\sqrt{C_*B}}{n^{Q/2}} .
\end{align*}
Let us put $Q=k+\frac{1}{2}$, then we obtain \eqref{4-1} with some constant $C^*.$

\bigskip

{\bf Case $p=2$ and $k=3$}. To illustrate \eqref{4-1} we consider the
expansion of the moments in the case $p=2$ and $k=3$. We can write
($\varepsilon =n^{-1/2}$)
\begin{align*}
n\Ex_{\vartheta_0}\Psi _{0,n}^2&=
\Ex_{\vartheta_0}\left(\psi_1\eta_n+\psi_2\eta_n^2\varepsilon
+\psi_3\eta_n^3\varepsilon^2\right)^2\\
&=\psi_1^2\Ex_{\vartheta_0}\eta_n^2+
2\psi_1\psi_2\Ex_{\vartheta_0}\eta_n^3\varepsilon+
\left[\psi_2^2+2\psi_1\psi_3\right]
\Ex_{\vartheta_0}\eta_n^4\varepsilon^2+q_n\varepsilon^3+p_n\varepsilon ^4,
\end{align*}
where $q_n $ and $p_n$ are bounded sequences. 
Remind that
\begin{align*}
\Ex_{\vartheta_0}\eta_n^2&=\int_0^Tg(t)^2\lambda(\vartheta_0,t)\,{\rm d}t
\equiv a_2,\\
\Ex_{\vartheta_0}\eta_n^3&=\frac{1}{\sqrt{n}}\int_0^Tg(t)^3\lambda(\vartheta_0,t)\,{\rm d}t
\equiv a_3\varepsilon,\\
\Ex_{\vartheta_0}\eta_n^4&=3\left(\int_0^Tg(t)^2
\lambda(\vartheta_0,t){\rm
  d}t\right)^2+\frac{1}{n}\int_0^Tg(t)^4\lambda(\vartheta_0,t)\,{\rm d}t\equiv
3a_2^2+a_4\varepsilon^2,
\end{align*} 
where we denoted 
\begin{align*}
a_m=\int_0^Tg(t)^m\lambda(\vartheta_0,t){\rm d}t,\qquad m=2,3,4.
\end{align*}
Therefore the first terms are 
\begin{align*}
n\Ex_{\vartheta_0}\Psi _{0,n}^2=
a_2\psi_1^2+\left[2a_3\psi_1\psi_2+3a_2^2\left(\psi_2^2+2\psi_1\psi_3\right)
\right]\varepsilon^2+\tilde q_n\varepsilon^3+\tilde p_n\varepsilon ^4,
\end{align*}
and we can write
\begin{align*}
\left|n\Ex_{\vartheta _0}\left(\check \vartheta _n-\vartheta
_0\right)^2-a_2\psi_1^2-\left[2a_3\psi_1\psi_2+3a_2^2\left(\psi_2^2+2\psi_1\psi_3\right) 
\right]\varepsilon^2\right|\leq C\varepsilon ^{5/2}
\end{align*}
or
\begin{align}
\label{4-3}
\frac{n}{\psi_1^2a_2}\Ex_{\vartheta _0}\left(\check \vartheta
_n-\vartheta_0\right)^2=1+\left[\frac{2\psi_2a_3}{\psi_1a_2}+\frac{3\psi_2^2a_2}{\psi
    _1^2}+\frac{6\psi_3a_2}{\psi _1} \right]
\frac{1}{n}+\frac{R_n}{n^{\frac{5}{4}}},
\end{align}
where $R_n$ is bounded sequence.
We see that the term of order $n^{-3/2}$ is absent.

{\bf Example 4}. Consider the model of $1$-periodic Poisson process $X_t,0\leq
t\leq T$ with intensity function
\begin{align*}
\lambda \left(\vartheta ,t\right)=2\sin\left(2\pi t+\vartheta \right)+3,\qquad
0\leq t\leq T=n, 
\end{align*}
where $\vartheta \in \left(\alpha ,\beta \right)$, $0<\alpha <\beta <\frac{\pi}{2}
$. Hence we have $n$ independent observations $X_j=\left(X_j\left(t\right),0\leq
t\leq 1\right)$, where  $X_j\left(t\right)=X_{j-1+t}-X_{j-1}$, $j=1,\ldots,n$.

Let us take $g\left(t\right)=\cos\left(2\pi t\right)$. Then 
\begin{align*}
m\left(\vartheta \right)& =\sin \vartheta ,\quad
\check\vartheta _n=\arcsin
\left(\frac{1}{n}\sum_{j=1}^{n}\int_{0}^{1}\cos\left(2\pi t\right) {\rm
  d}X_j\left(t\right) \right),\\
G\left(y\right)&=\arcsin y,\qquad a_2=\frac{3}{2},\qquad a_3=\frac{3}{4}\sin \vartheta ,\\
\psi _1&=\frac{1}{\cos \vartheta },\qquad \psi _2=\frac{\sin \vartheta
}{2\cos^3 \vartheta },\qquad  \psi _3=\frac{1+2\sin^2 \vartheta}{6\cos^5
  \vartheta }.
\end{align*}
Suppose that the true value is $\vartheta _0=\frac{\pi }{3}$. Then
\begin{align*}
a_2=\frac{3}{2},\qquad a_3=\frac{3\sqrt{3}}{8},\qquad \psi _1=2,\qquad \psi
_2=2\sqrt{3},\qquad  \psi _3=\frac{40}{3}.
\end{align*}
The expansion of the second moment is
\begin{align*}
n\Ex_{\frac{\pi }{3}}\left(\check\vartheta _n-\frac{\pi
}{3}\right)^2=6+\frac{450}{n}+\frac{R_n}{n^{\frac{5}{4}}} 
\end{align*}
{\bf Simulation.} Let us check the last expansion by numerical
simulations. The limit value of the right-hand side is 6. Suppose that
$n=1000$, then
\begin{align*}
1000\,\Ex_{\frac{\pi }{3}}\left(\check\vartheta _{1000}-\frac{\pi
}{3}\right)^2=6+\frac{450}{1000}+\frac{R_{1000}}{10^{\frac{15}{4}}} \approx 6,5.
\end{align*}
The estimator $\check\vartheta _{1000} $ was simulated 10 000 times,  $\check\vartheta _{m,1000},m=1,\ldots, 10000
$. The empirical error obtained by this simulation 
\begin{align*}
\frac{1}{N}\sum_{m=1}^{N}\left(\check\vartheta _{m,1000}-\frac{\pi
}{3}\right)^2=6,53
\end{align*}
corresponds well to the value obtained by asymptotic expansion.

\section{Expansion of the distribution function}

We discuss below the non asymptotic expansions of the distribution function
of the MME $\check\vartheta _n$. We consider the cases $k=1$ and  $k=2$
only, which correspond to the representations
\begin{align}
\label{5-1}
\check\vartheta _n&=\vartheta
_0+\left(\psi_1\eta_n\varepsilon+\psi_2\eta_n^2\varepsilon^2+\bar r_n
\varepsilon^{5/2}\right)
\1_{\left\{\AA\right\}}+\phi_n\1_{\left\{\AA^c\right\}} ,\\ \check\vartheta
_n&=\vartheta
_0+\left(\psi_1\eta_n\varepsilon+\psi_2\eta_n^2\varepsilon^2+\psi_3\eta_n^3\varepsilon^3+\bar 
r_n \varepsilon^{7/2}\right) 
\1_{\left\{\AA\right\}}+\phi_n\1_{\left\{\AA^c\right\}},
\label{5-2}
\end{align}
where $\varepsilon =n^{-1/2}$, $\left|\bar r_n\right|<1,
\left|\phi_n\right|<\beta -\alpha $ and we have estimates  \eqref{3-2}. Of course, $ \psi
_l=\psi _l\left(\vartheta _0\right),l=1,2,3 $, $\eta _n=\eta _n\left(\vartheta
_0\right)$ and $\bar r _n=\bar r_n\left(\vartheta _0\right)$  but we omit this
dependence to simplify the exposition. Moreover we keep the same notation for
the different random variables $\bar r_n$,$\phi_n$  and  the sets $\AA$ in
 \eqref{5-1} and \eqref{5-2}. Our goal is to describe the first non
Gaussian terms of these expansions. The proof in the case $k=2$ we give in
details. The proof in the  case $k=3$ is much more cumbersome that is why we
  obtain just formally the first two terms after Gaussian and for detailed proofs
  propose to read the section 3.4 in \cite{Kut98}, where the similar
  expansions were studied. 

 Introduce the notation
\begin{align*}
\xi _n&=\frac{1}{\sqrt{a_2n}}\sum_{j=1}^{n}\int_{\cal T}^{}g\left(t\right){\rm
  d}\pi _j\left(t\right)=\frac{\eta _n }{\sqrt{a_2}}\Longrightarrow {\cal
  N}\left(0,1\right),\\
b_2&=\frac{\psi _2\sqrt{a_2}}{\psi _1},\qquad b_3=\frac{\psi _3a_2}{\psi
  _1},\qquad \tilde r_n=\frac{\bar r_n}{\psi _1\sqrt{a_2}},\qquad \tilde
\phi_n=\frac{\phi_n\sqrt{n}}{\psi _1\sqrt{a_2}}. 
\end{align*}
Then the stochastic expansion can be written as follows
\begin{align*}
\sqrt{\frac{n}{\psi _1^2a_2}}\left(\check\vartheta _n-\vartheta _0\right)=\left\{\xi
_n+b_2\xi _n^2\varepsilon +b_3\xi _n^3\varepsilon ^2 +\tilde r_n \varepsilon ^{5/2}\right\}\1_{\left\{\AA\right\}}+\tilde\phi_n\1_{\left\{\AA^c\right\}}  
\end{align*}
Our goal is to obtain the expansion of the probability
\begin{align*}
F_n\left(x\right)&=\Pb_{\vartheta _0}\left( \sqrt{\frac{n}{\psi _1^2a_2}}\left(\check\vartheta
_n-\vartheta _0\right)<x\right)\\
& =\Pb_{\vartheta _0}\left(\xi
_n+b_2\xi _n^2\varepsilon +b_3\xi _n^3\varepsilon ^2 +\tilde r_n \varepsilon
^{5/2}  <x,\AA\right) + \Pb_{\vartheta _0}(\tilde\phi_n<x,\AA^c)
\end{align*}
by the powers of $\varepsilon $. Let us denote 
\begin{align*}
\Phi _{n,\varepsilon }=\xi
_n+b_2\xi _n^2\varepsilon +b_3\xi _n^3\varepsilon ^2 +\tilde r_n \varepsilon^{5/2},\qquad \Phi _{0,\varepsilon }= \xi _n+b_2\xi _n^2\varepsilon
  +b_3\xi _n^3\varepsilon ^2.
\end{align*}
Then we can write
\begin{align*}
\Pb_{\vartheta _0}\left(\Phi _{n,\varepsilon }<x,\AA\right)\leq    F_n\left(x\right)\leq \Pb_{\vartheta _0}\left(\Phi _{n,\varepsilon }<x,\AA\right)+\Pb_{\vartheta _0}(\AA^c).
\end{align*}
Therefore, it is sufficient to study the probability $\Pb_{\vartheta
  _0}\left(\Phi _{n,\varepsilon }<x,\AA\right) $. Recall that  on the set  $\AA$
we have 
$\left|\tilde r_n\right|\leq K=\left(\psi _1\sqrt{a_2}\right)^{-1}$. Hence
\begin{align*}
\Pb_{\vartheta _0}\left(\Phi _{n,\varepsilon }<x,\AA\right)&\leq \Pb_{\vartheta
  _0}\left( \xi _n+b_2\xi _n^2\varepsilon +b_3\xi _n^3\varepsilon ^2
<x+K\varepsilon ^{5/2},\AA\right)\nonumber\\
&=\Pb_{\vartheta _0}\left(\Phi _{0,\varepsilon }<x+K\varepsilon ^{5/2},\AA\right)
\end{align*} 
and
\begin{align*}
\Pb_{\vartheta _0}\left(\Phi _{n,\varepsilon }<x,\AA\right)&\geq \Pb_{\vartheta
  _0}\left( \xi _n+b_2\xi _n^2\varepsilon +b_3\xi _n^3\varepsilon ^2
<x-K\varepsilon ^{5/2},\AA\right)\nonumber\\
&=\Pb_{\vartheta _0}\left(\Phi _{0,\varepsilon }<x-K\varepsilon ^{5/2},\AA\right).
\end{align*}
Further, we can write 
\begin{align*}
&\Pb_{\vartheta _0}\left( \Phi _{0,\varepsilon } <y,\AA\right)\leq
\Pb_{\vartheta _0}\left(\Phi _{0,\varepsilon } <y\right),\\ 
&\Pb_{\vartheta
  _0}\left( \Phi _{0,\varepsilon } <y,\AA\right)\geq \Pb_{\vartheta _0}\left(
\Phi _{0,\varepsilon } <y\right)-\Pb_{\vartheta _0}\left(\AA^c\right).
\end{align*}

Let us study the probabilities
\begin{align*}
P_n^{\left(1\right)}\left(y\right)&=\Pb_{\vartheta _0}\left(\xi _n+b_2\xi
_n^2\varepsilon <y\right),\\
P_n^{\left(2\right)}\left(y\right)&=\Pb_{\vartheta _0}\left(\xi _n+b_2\xi
_n^2\varepsilon+b_3\xi _n^3\varepsilon ^2 <y\right)
\end{align*}
separately. 

  The expansion of the probability
$P_n^{\left(1\right)}\left(\cdot \right) $ we give with proof and for the
probability $P_n^{\left(2\right)}\left(\cdot \right) $ we present a formal expansion
without detailed description of the reminder term.

\subsection{ Expansion of $P_n^{\left(1\right)}\left(y\right)$}

\begin{proposition}
\label{P5-1} Suppose that   the conditions ${\cal L}$, $\psi _1>0$ and 
\begin{align}
\label{5-3}
\inf_{a_2b\left|v\right|>1}\int_{0}^{\tau
}\sin^2\left(vg\left(t\right)\right)\lambda \left(\vartheta _0,t\right){\rm
  d}t >0
\end{align}
hold. Then we have the expansion
\begin{align}
\label{5-4}
&\Pb_{\vartheta _0}\left(\frac{\sqrt{n}\left(\check\vartheta _n-\vartheta _0
  \right)}{\psi _1\sqrt{a_2}}<y\right)\nonumber\\ &\qquad \qquad \qquad \qquad
=F\left(y\right)+\frac{\left[a_3- \left(a_3+8a_2^{{3/2}}\,b_2\right)y^2\right]}{6
  a_2^{3/2}} f\left(y\right)\varepsilon +o\left(\varepsilon \right) .
\end{align}

\end{proposition}
{\bf Proof.}  First we recall the expansion of the distribution function of
the stochastic integral (see \cite{Kut98}, Theorem 3.4). We
use the same notation as above
\begin{align*}
\xi _n&=\frac{1}{\sqrt{a_2n}}\sum_{j=1}^{n}\int_{0}^{\tau
}g\left(t\right){\rm d}\pi _j\left(t\right) ,\qquad b_2=\frac{\psi _2\sqrt{a_2}}{\psi _1},\quad \varepsilon =\frac{1}{\sqrt{n}},\\
 a_m&=\int_{0}^{\tau
}g\left(t\right)^m\lambda \left(\vartheta
  _0,t\right){\rm d}t,\quad m=2,3,4.
\end{align*}

 Suppose that the function
  $g\left(t\right)$   is non constant and  satisfies  \eqref{3-0}.
Then
\begin{align*}
\Pb_{\vartheta _0}\left(\xi _n <x\right)=F\left(x\right)+\frac{a_3}{6
  a_2^{3/2}\sqrt{2\pi }} \left(1-x^2\right)e^{-x^2/2}\varepsilon
+o\left(\varepsilon \right).
\end{align*}
We have
\begin{align*}
P_n^{\left(1\right)}\left(y\right)=\int_{x+b_2x^2\varepsilon <y}^{}{\rm
  d}\Pb_{\vartheta _0}\left(\xi _n <x\right) .
\end{align*}
The condition $x+b_2x^2\varepsilon <y $ is equivalent to $b_2\varepsilon
\left(x-x_1\right)\left(x-x_2\right)<0 $, where the roots are 
\begin{align*}
x_1=\frac{-1-\sqrt{1+4b_2\varepsilon y}}{2b_2\varepsilon },\qquad
x_2=\frac{-1+\sqrt{1+4b_2\varepsilon y}}{2b_2\varepsilon }. 
\end{align*}
Suppose that $b_2>0$ and $y\in \YY$, where $\YY$ is some compact. Then for
sufficiently small $\varepsilon $ we have 
\begin{align*}
x_1=-\frac{1}{b_2\varepsilon }\left(1+O\left(\varepsilon
^2\right)\right),\qquad x_2= y-b_2y^2\varepsilon +O\left(\varepsilon
^2\right),
\end{align*}
and
\begin{align*}
P_n^{\left(1\right)}\left(y\right)&=\Pb_{\vartheta _0}\left(\xi _n
<y-b_2y ^2\varepsilon +O\left(\varepsilon ^2\right)\right)\\
&\qquad \qquad -\Pb_{\vartheta _0}\left(\xi _n
<- \left(b_2\varepsilon\right)^{-1} \left(1+O\left(\varepsilon
^2\right)\right)\right)\\
&=\Pb_{\vartheta _0}\left(\xi _n
<y-b_2y ^2\varepsilon \right)+O\left(\varepsilon ^2\right)
\end{align*}
because   $\Pb_{\vartheta _0}\left(\xi _n
<- \left(b_2\varepsilon\right)^{-1} \left(1+O\left(\varepsilon
^2\right)\right)\right)$ is exponentially small. 

Further
\begin{align*}
F\left(y-b_2y ^2\varepsilon
\right)=F\left(y\right)-b_2y ^2f\left(y\right)\varepsilon
+O\left(\varepsilon ^2\right) .
\end{align*}
Hence for $b_2>0$ we have
\begin{align}
\label{5-5}
P_n^{\left(1\right)}\left(y\right)&=F\left(y\right)+\left[\frac{a_3}{6
  a_2^{3/2}} \left(1-y^2\right)-b_2y ^2 \right]   f\left(y\right)\varepsilon
+o\left(\varepsilon \right)  .
\end{align}
If $b_2<0$, then the condition $x+b_2x^2\varepsilon <y$ is equivalent to two
conditions 
\begin{align*}
x<x_2=y-b_2y^2\varepsilon +O\left(\varepsilon ^2\right)\; {\rm  and}\quad 
x>x_1=-\left( b_2\varepsilon \right)^{-1}\left(1+O\left(\varepsilon
^2\right)\right)\rightarrow +\infty .
\end{align*}
Therefore
\begin{align*}
P_n^{\left(1\right)}\left(y\right)&=\Pb_{\vartheta _0}\left(\xi _n
<y-b_2y ^2\varepsilon +O\left(\varepsilon ^2\right)\right)\\
&\qquad \qquad +\Pb_{\vartheta _0}\left(\xi _n
>- \left(b_2\varepsilon\right)^{-1} \left(1+O\left(\varepsilon
^2\right)\right)\right)\\
&=\Pb_{\vartheta _0}\left(\xi _n
<y-b_2y ^2\varepsilon \right)+O\left(\varepsilon ^2\right)
\end{align*}
and we obtain the same expansion \eqref{5-5}.

The expansion of the corresponding density function  we
obtain by formal derivation of the first two terms of the right hand side of
\eqref{5-4} 
\begin{align}
\label{5-6}
f_{0,\varepsilon }^{\left(1\right)}\left(y\right)=f\left(y\right)+\left[b_2
y+B_3\left(y^3-3y\right)\right]f\left(y\right) \varepsilon +o\left(\varepsilon \right),
\end{align}
where we denoted $B_3=a_3a_2^{-3/2}6^{-1}+b_2$.

\subsection{ Expansion of $P_n^{\left(2\right)}\left(y\right)$}

Now we consider the probability
\begin{align*}
\Pb_{\vartheta _0}\left(\frac{\sqrt{n}\left(\check\vartheta
_n-\vartheta _0 \right)}{\psi _1\sqrt{a_2}}<y\right)
=\Pb_{\vartheta _0}\left(\Phi _{0,\varepsilon }<y\right)+o\left(\varepsilon ^2\right) ,
\end{align*}
where  $\Pb_{\vartheta _0}\left(\Phi _{0,\varepsilon}<y\right)=\Pb_{\vartheta _0}\left(\xi _n+b_2\xi
_n^2\varepsilon +b_3\xi _n^3 \varepsilon ^2<y\right)=P_n^{\left(2\right)}\left(y\right)  $.

 The  characteristic function $ M_\varepsilon
  (v)=\Ex_{\vartheta _0}\exp({{\rm i} v\Phi _{0,\varepsilon }})$  can be
  expanded by the powers of $\varepsilon $ as follows
\begin{align}
\label{5-7}
M_\varepsilon (v)&=\Ex_{\vartheta _0}e^{{\rm i} v\xi  _{n }}\left(1+{\rm i}vb_2\xi
_n^2\varepsilon +{\rm i}vb_3\xi _n^3\varepsilon ^2 +\frac{\left({\rm i} v\right)^2}{2}b_2^2\xi
_n^4\varepsilon ^2+h\left(v,\xi _n\right)\varepsilon ^3\right)\nonumber \\
&=\Ex_{\vartheta _0}e^{{\rm i} v\xi  _{n }}+{\rm i}vb_2 \Ex_{\vartheta _0}\xi  _{n
}^2e^{{\rm i} v\xi  _{n }}\,\varepsilon +{\rm i}vb_3 \Ex_{\vartheta _0}\xi  _{n
}^3e^{{\rm i} v\xi  _{n }}\,\varepsilon ^2\nonumber \\
&\qquad \qquad  +\frac{\left({\rm i}v\right)^2}{2}b_2^2 \Ex_{\vartheta _0}   \xi
_n^4e^{{\rm i} v\xi  _{n }}\varepsilon ^2+H_n\left(v\right)\varepsilon ^3.
\end{align}
Here $h\left(\cdot \right)$ and $H_n\left(\cdot \right)$ are the corresponding
resiudals.

By direct calculation  we obtain the values (below
$g_t=a_2^{-1/2}g\left(t\right)$ and $\lambda _t=\lambda \left(\vartheta
_0,t\right)$)
\begin{align*}
&\Ex_{\vartheta _0}e^{{\rm i} v\xi _{n }}=\exp\left\{\varepsilon
^{-2}\int_{\cal T}^{}\left(e^{{\rm i} vg_t\varepsilon }-1-{\rm
  i} vg_t\varepsilon \right)  \lambda _t   {\rm d}t\right\}\equiv {\cal E}\left(v,\varepsilon 
\right),\\ 
&\Ex_{\vartheta _0}\xi _ne^{{\rm i} v\xi _{n }}=\left[\varepsilon
  ^{-1}\int_{\cal T}^{}\left(e^{{\rm i} vg_t\varepsilon
  }-1\right)g_t\lambda _t{\rm
    d}t\right]{\cal E}\left(v,\varepsilon \right),\\ 
&\Ex_{\vartheta _0}\xi _n^2e^{{\rm i} v\xi _{n }}=\left[J_\varepsilon
  \left(v\right)^2 +\int_{\cal T}^{}e^{{\rm i} vg_t\varepsilon 
  }g_t^2\lambda _t{\rm d}t\right] {\cal E}\left(v,\varepsilon \right),
\end{align*}
where we denoted
\begin{align*}
J_\varepsilon \left(v\right)=\varepsilon
  ^{-1}\int_{\cal T}^{}\left(e^{{\rm i} vg_t\varepsilon
  }-1\right)g_t\lambda _t{\rm
    d}t.
\end{align*}
The similar but cumbersome expressions we have for expectations $\Ex_{\vartheta _0}\xi
_n^3e^{{\rm i} v\xi _{n }} $ and $\Ex_{\vartheta _0}\xi
_n^4e^{{\rm i} v\xi _{n }} $. As we need just the first terms corresponding
$\varepsilon ^0$ we do not give here these expressions. 

According to the expression \eqref{5-7} we need expansion of $\Ex_{\vartheta
  _0}e^{{\rm i} v\xi _{n }}$ up to $\varepsilon ^2$, the expansion of
$\Ex_{\vartheta _0}\xi _n^2e^{{\rm i} v\xi _{n }} $ up to $\varepsilon $ and
expansions of $\Ex_{\vartheta _0}\xi _n^3e^{{\rm i} v\xi _{n }} $ and
$\Ex_{\vartheta _0}\xi _n^4e^{{\rm i} v\xi _{n }} $ up to the first term
$\varepsilon ^0$ only.  Let us denote
\begin{align*}
\hat a_m=a_2^{-m/2}\int_{\cal T}^{}g\left(t\right)^m\lambda \left(\vartheta
_0,t\right){\rm d}t ,\qquad m=3,4,\qquad \left(\hat a_2=1\right)
\end{align*}
and remark, that $J_\varepsilon \left(v\right)={\rm i} v+\frac{\left({\rm i}
  v\right)^2}{2} \hat a_3\varepsilon +O\left(\varepsilon ^2\right) $.

We have
\begin{align*}
 \Ex_{\vartheta _0}e^{{\rm i} v\xi _{n
 }}&=e^{-\frac{v^2}{2}}\left(1+\frac{\left({\rm i}v\right)^3 }{6}\hat
 a_3\varepsilon +\frac{\left({\rm i}v\right)^4}{24}\hat a_4\varepsilon
 ^2+\frac{\left({\rm i}v\right)^6}{72}\hat a_3^2\varepsilon
 ^2+O\left(\varepsilon ^3\right) \right),\\
 \Ex_{\vartheta _0}\xi _n^2e^{{\rm
     i} v\xi _{n }}&=e^{-\frac{v^2}{2}}\left(1+\left({\rm
   i}v\right)^2+\left[\left({\rm i}v\right)+\frac{7\left({\rm
       i}v\right)^3}{6}+\frac{\left({\rm i}v\right)^5}{6} \right] \hat a_3
 \varepsilon\right) +O\left(\varepsilon ^2\right),\\ 
\Ex_{\vartheta
   _0}\xi _n^3e^{{\rm i} v\xi _{n }}&=e^{-\frac{v^2}{2}}\left[\left({\rm i}
  v\right)^3+3\left({\rm i}     v\right)\right]+O\left(\varepsilon \right),\\ 
\Ex_{\vartheta _0}\xi 
 _n^4e^{{\rm i} v\xi _{n }}&=e^{-\frac{v^2}{2}}\left[\left({\rm i} v\right)^4+6\left({\rm i}
    v\right)^2+3\right]+O\left(\varepsilon \right). 
\end{align*}

These expressions allow us to write the expansion of the characteristic
function \eqref{5-7}
\begin{align*}
M_\varepsilon (v)&=e^{-\frac{v^2}{2}}+ \left({\rm
  i}v\right)e^{-\frac{v^2}{2}}B_1 \varepsilon+
\left({\rm i}v\right)^2e^{-\frac{v^2}{2}}B_2\varepsilon^2 +\left({\rm
  i}v\right)^3e^{-\frac{v^2}{2}}B_3 \varepsilon\\
&\qquad +\left({\rm i}v\right)^4e^{-\frac{v^2}{2}}B_4
\varepsilon^2 +\left({\rm
  i}v\right)^6e^{-\frac{v^2}{2}}B_6\varepsilon ^2 +R_\varepsilon
\left(v\right)\varepsilon ^3,
\end{align*}
where $R_\varepsilon \left(v\right)\varepsilon ^3$ is the corresponding
resiudal and we used the notation
\begin{align*}
B_1&=b_2, \qquad B_2=b_2\hat a_3+\frac{3}{2}b_2^2+3b_3,\qquad B_3=\frac{1
}{6}\hat  a_3 +b_2,\\
B_4&=\frac{1}{24}\hat a_4+ \frac{7}{6}b_2\hat a_3+b_3+3b_2^2,\qquad 
B_6=\frac{1}{72}\hat  a_3^2+\frac{1}{6}b_2\hat a_3+\frac{1}{2}b_2^2. 
\end{align*}

Recall the relation
\begin{align}
\label{5-8}
\frac{1}{2\pi }\int_{\cal R}^{}e^{-{\rm i}vx}\left({\rm
  i}v\right)^me^{-\frac{v^2}{2}}{\rm d}v=H_m\left(x\right)f\left(x\right) 
,\quad  m=1,2,\ldots
\end{align}
where $H_m\left(\cdot \right)$ are Hermite polynomials. We have
\begin{align}
\label{5-9}
  H_1(x)&=x,\qquad H_2(x)=x^2-1,\qquad H_3(x)=x^3-3x,\nonumber\\
H_4(x)&=x^4-6x^2+3, \qquad H_6(x)=x^6-15x^4+45x^2-15 .
\end{align}
Let us denote
\begin{align*}
M_{0,\varepsilon }(v)&=e^{-\frac{v^2}{2}}+ \left({\rm
  i}v\right)e^{-\frac{v^2}{2}}B_1 \varepsilon+
\left({\rm i}v\right)^2e^{-\frac{v^2}{2}}B_2\varepsilon^2 +\left({\rm
  i}v\right)^3e^{-\frac{v^2}{2}}B_3 \varepsilon\\
&\qquad +\left({\rm i}v\right)^4e^{-\frac{v^2}{2}}B_4
\varepsilon^2 +\left({\rm
  i}v\right)^6e^{-\frac{v^2}{2}}B_6\varepsilon ^2.
\end{align*}
Using these relations we obtain the inverse Fourier transform of
$M_{0,\varepsilon} (v)$
\begin{align*}
f_{0,\varepsilon }\left(x\right)&=\frac{1}{2\pi }\int_{\cal R}^{}e^{-{\rm
    i}vx}M_{0,\varepsilon} (v){\rm d}v
=f\left(x\right)+B_1H_1\left(x\right)f\left(x\right)\varepsilon \\
&\quad  +B_3H_3\left(x\right)f\left(x\right)\varepsilon +
\left[B_2H_2\left(x\right)+B_4H_4\left(x\right)+B_6H_6\left(x\right)
  \right]f\left(x\right) \varepsilon ^2 .
\end{align*}
Introduce the function 
\begin{align*}
F_{0,\varepsilon }\left(x\right)&=\int_{-\infty }^{x}f_{0,\varepsilon
}\left(y\right){\rm d}y .
\end{align*}
It is possible to show that 
\begin{align*}
F_{n,\varepsilon }\left(x\right)=\Pb_{\vartheta _0}\left(\sqrt{\frac{n}{\psi _1^2a_2}}\left(\check\vartheta
_n-\vartheta _0\right)<x\right)=F_{0,\varepsilon 
}\left(x\right)+ O\left(\varepsilon ^{\frac{5}{2}}\right).
\end{align*}
 The proof of the convergence
$O\left(\varepsilon ^{\frac{5}{2}}\right) $ follows the same steps as the proof
of the Theorem 3.4 in \cite{Kut98}.
 
Therefore the functions 
\begin{align}
\label{5-10}
F_{0,\varepsilon }^{\left(1\right)}\left(x\right)&=\frac{1}{\sqrt{2\pi }}\int_{-\infty
}^{x}e^{-\frac{y^2}{y}}{\rm d}y +\int_{-\infty
}^{x}\left[B_1H_1\left(y\right)+B_3  H_3\left(y\right)\right]
f\left(y\right){\rm d}y \,\frac{1}{\sqrt{n}}  
\end{align}
and
 \begin{align*}
F_{0,\varepsilon }^{\left(2\right)}\left(x\right)&=\frac{1}{\sqrt{2\pi }}\int_{-\infty
}^{x}e^{-\frac{y^2}{y}}{\rm d}y +\int_{-\infty
}^{x} \left[B_1H_1\left(y\right)+B_3  H_3\left(y\right)\right]
f\left(y\right){\rm d}y \,\frac{1}{\sqrt{n}} \\ 
&\quad +\int_{-\infty
}^{x}\left[B_2H_2\left(y\right)+B_4H_4\left(y\right)+B_6H_6\left(y\right)
  \right] f\left(y\right){\rm d}y \,\frac{1}{n}
\end{align*}
can be considered as approximations of the distribution function of the MME
with the corresponding precisions respectively. Note that the expression
\eqref{5-10} corresponds well to the expansion \eqref{5-6} obtained
before. 

It is interesting to obtain the approximation of the second moment of the MME
\begin{align*}
\frac{n}{\psi _1^2a_2}\Ex_{\vartheta _0}\left(\check\vartheta _n-\vartheta
_0\right)^2=\int_{ \cal R}^{}x^2{\rm d}F_{0,\varepsilon
}\left(x\right)+O\left(\varepsilon ^{\frac{5}{2}}\right) 
\end{align*}
with the help of the approximation of the distribution function. We have
\begin{align*}
\int_{\cal R}^{}x^2f_{0,\varepsilon }\left(x\right){\rm d}x&= 1+B_2\int_{\cal
  R}^{}x^2H_2\left(x\right)f\left(x\right){\rm d}x\, \frac{1}{n}= 1+\frac{2B_2}{n}\\
&=1+\left(2b_2\hat a_3+3b_2^2+6b_3\right)\frac{1}{n}\\
&=1+\left(\frac{2\psi _2 a_3}{\psi _1a_2}+\frac{3\psi _2^2a_2}{\psi
  _1^2}+\frac{6\psi _3a_2}{\psi _1}\right)\frac{1}{n} .
\end{align*}
Here we used the properties of the Hermite polynomials
\begin{align*}
\int_{\cal R}^{}x^2H_m\left(x\right)p_1\left(x\right){\rm d}x=0,\qquad m=1,3,4,6.
\end{align*}
Hence
\begin{align}
\frac{n}{\psi _1^2a_2}\Ex_{\vartheta _0}\left(\check\vartheta _n-\vartheta
_0\right)^2=1+\left(\frac{2\psi _2 a_3}{\psi _1a_2}+\frac{3\psi _2^2a_2}{\psi 
  _1^2}+\frac{6\psi _3a_2}{\psi _1}\right)\frac{1}{n} +O\left(n
^{-\frac{5}{4}}\right) . 
\label{5-11}
\end{align}
We see that the expressions \eqref{4-3} and \eqref{5-11} coincide.

{\bf Example 4} Suppose that the observed Poisson process is from the Example
4 and the true value is $\vartheta _0=\frac{\pi }{3}$. Then 
\begin{align*}
B_1\left(\frac{\pi }{3}\right)=\frac{3\sqrt{2}}{2},\qquad \quad
B_3\left(\frac{\pi }{3}\right) =\frac{11\sqrt{2}}{6}
\end{align*}
and we can write the approximation
\begin{align*}
&\Pb_{\vartheta _0}\left(\sqrt{\frac{2n\cos \vartheta
      _0}{3}}\left(\check\vartheta 
_n-\vartheta _0\right)<x\right)\\
&\qquad \quad \approx \frac{1}{\sqrt{2\pi }}\int_{-\infty
}^{x}e^{-\frac{y^2}{2}}{\rm d}y+\frac{B_1\left(\vartheta
  _0\right)}{\sqrt{2\pi }}\int_{-\infty 
}^{x}ye^{-\frac{y^2}{2}}{\rm d}y\;\frac{1}{\sqrt{n}}\\
&\qquad\qquad\quad \qquad +\frac{B_3\left(\vartheta
  _0\right)}{\sqrt{2\pi }}\int_{-\infty 
}^{x}\left(y^3-3y\right)e^{-\frac{y^2}{2}}{\rm d}y\;\frac{1}{\sqrt{n}}\\
&\qquad\quad = \frac{1}{\sqrt{2\pi }}\int_{-\infty
}^{x}e^{-\frac{y^2}{2}}{\rm d}y+\frac{3}{2\sqrt{\pi }}\int_{-\infty 
}^{x}ye^{-\frac{y^2}{2}}{\rm d}y\;\frac{1}{\sqrt{n}}\\
&\qquad\qquad\quad \qquad +\frac{11}{6\sqrt{\pi }}\int_{-\infty 
}^{x}\left(y^3-3y\right)e^{-\frac{y^2}{2}}{\rm d}y\;\frac{1}{\sqrt{n}}.
\end{align*}
We have no approximation of the density function but can write formally the
 density of the approximating distribution function
\begin{align*}
f_{0,n}\left(x\right)&=\frac{1}{\sqrt{2\pi
}}e^{-\frac{x^2}{2}}+\left[\frac{3 x}{2\sqrt{\pi
}}e^{-\frac{x^2}{2}}+\frac{11\left(x^3-3x\right)}{6\sqrt{\pi
}}e^{-\frac{x^2}{2}}\right]\;\frac{1}{\sqrt{n}} .
\end{align*}

 \begin{figure}[ht]
  \hspace{4cm}\includegraphics[width=7cm]   {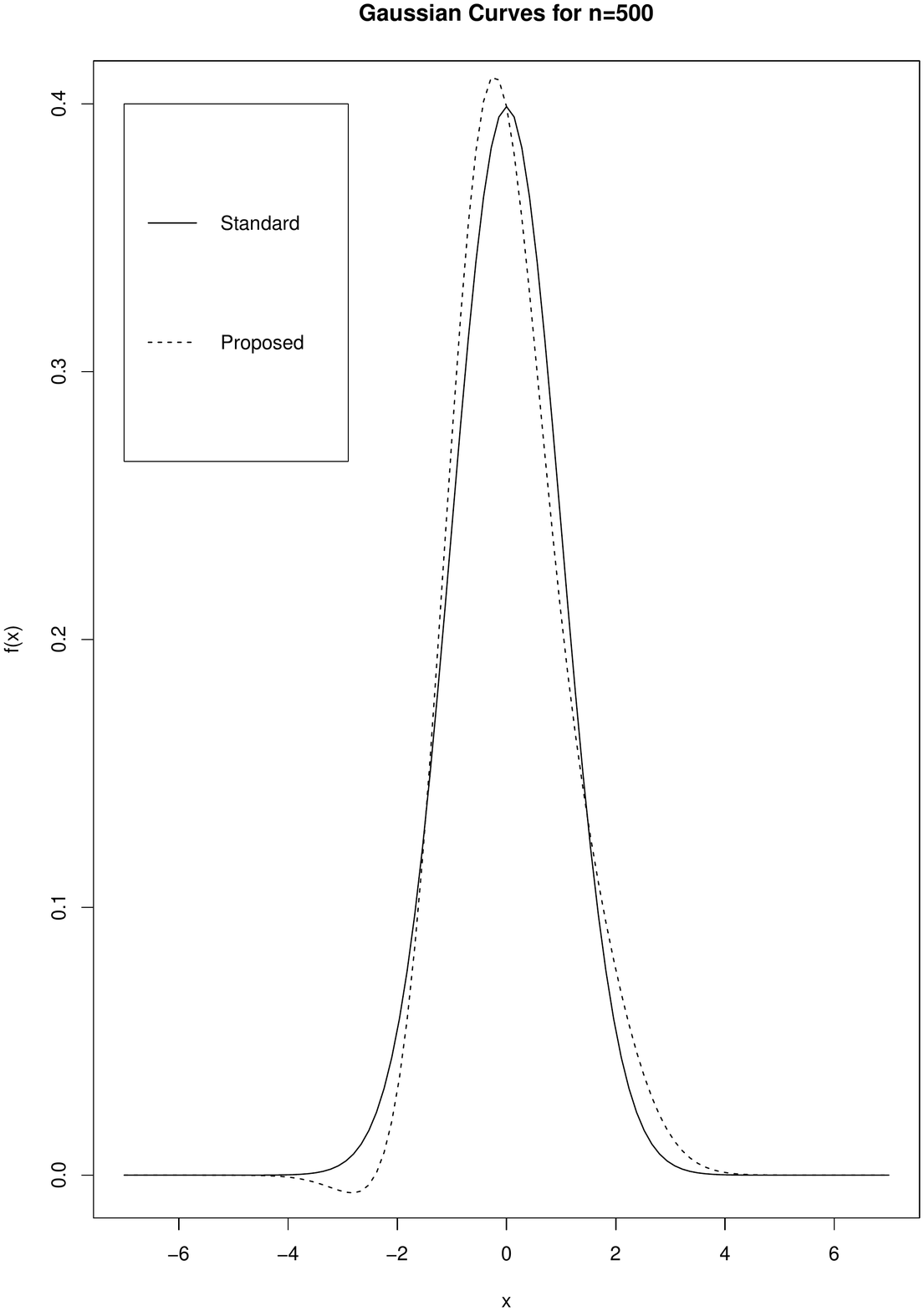}
  \caption{\it Densities: standard Gaussian $f\left(\cdot \right)$ and  proposed
    approximation $f_{0,500}\left(\cdot \right)$.} 
  \end{figure}

\section{Discussion}

The presented here results admit several generalizations. The case of vector
parameter $\vartheta $ can be treated by a similar way, but the exposition
will be more complicate. The main idea - to have non asymptotic description of
estimators can be applied to many different statistical models and many
statistical estimators. It was already realized for some models of continuous
time processes (Gaussian, diffusion, inhomogeneous Poisson) but we suppose
that many other models, for example, time series can be studied using the same
ideas.

\section*{Acknowledgments} This research was financially supported by the
Ministry of Education and Science of the Russian Federation (research project
No. FSWF-2020-0022).

\end{document}